\documentclass[10pt,a4paper]{article}

\usepackage[a4paper,textwidth=380pt,textheight=569.68pt]{geometry}
\usepackage{bbm}
\usepackage{amsfonts}
\usepackage{amsmath}
\usepackage{amssymb}
\usepackage{amsthm}
\usepackage{newtxtext}
\usepackage{newtxmath}
\DeclareMathAlphabet{\mathcal}{OMS}{cmsy}{m}{n}
\DeclareSymbolFont{CMlargesymbols}{OMX}{cmex}{m}{n}
\DeclareMathDelimiter{(}{\mathopen} {operators}{"28}{CMlargesymbols}{"00}
\DeclareMathDelimiter{)}{\mathclose}{operators}{"29}{CMlargesymbols}{"01}
\usepackage{enumitem}
\usepackage{graphicx}
\usepackage{natbib}
\usepackage[pdfborder={0 0 0}]{hyperref}

\theoremstyle{definition}
\newtheorem{lemma}{Lemma}
\newtheorem{example}{Example}
\newtheorem{definition}{Definition}

\def\e{\mathrm{e}}
\def\E{E}
\DeclareMathOperator{\Leb}{Leb}
\def\N{\mathbb{N}}
\DeclareMathOperator{\Poiss}{Poiss}
\def\R{\mathbb{R}}

\newcommand{\rmd}{\mathrm{d}}
\newcommand{\PP}{P}
\DeclareMathOperator{\Var}{Var}
\newcommand{\myp}{\mbox{$\:\!$}}

\begin{document}
\bibliographystyle{jochen}

\title{Exact Sampling of Almost Surely Finite Sequences of Independent Events}

\author{%
  Leonid V.\ Bogachev\thanks{\,\href{mailto:L.V.Bogachev@leeds.ac.uk}{L.V.Bogachev@leeds.ac.uk}}\\[.3pc]
  \normalsize\it School of Mathematics\\
  \normalsize\it University of Leeds
  \and
  Ruheyan Nuermaimaiti\thanks{\,\href{mailto:r.nuermaimaiti@imperial.ac.uk}{r.nuermaimaiti@imperial.ac.uk}}\\[.3pc]
  \normalsize\it Department of Mathematics\\
  \normalsize\it Imperial College London
  \and
  Jochen Voss\thanks{\,Corresponding author, \href{mailto:J.Voss@leeds.ac.uk}{J.Voss@leeds.ac.uk}}\\[.3pc]
  \normalsize\it School of Mathematics\\
  \normalsize\it University of Leeds
}
\date{November 26, 2026}
\maketitle

\begin{abstract}
  We consider infinite sequences of independent events whose probabilities are
  summable, so that almost surely only finitely many of the events occur.
  We present an algorithm which exactly samples the set of events which
  occur, using a finite number of steps.  The algorithm works by
  generating proposals from a Poisson process, followed by a rejection step.
  We illustrate use of the algorithm with the help of
  three applications: prime factorisation of
  zeta-distributed random integers, random integer partitions and convex
  lattice polygonal lines.


\end{abstract}

\section{Introduction}

Consider a sequence of independent Bernoulli random variables $(X_n)_{n\in\N}$,
such that we have $\PP(X_n = 1) = 1 - \PP(X_n = 0) = p_n$ for all $n\in\N$.
If we assume that $\sum_{n=1}^\infty p_n < \infty$, the Borel--Cantelli lemma guarantees that, almost surely (a.s.), the events $\{X_n = 1\}$ happen only for finitely many~$n$.
Thus, while there are infinitely many $X_n$ in the sequence, an exact
representation of a random sample from the sequence can be stored using
only a finite amount of computer memory.  Here we show
that in many cases it is possible to exactly sample random instances of
$(X_n)_{n\in\N}$.  The main obstacle to overcome is the fact that
the position of the last non-zero $X_n$ is determined by the infinite tail of the random sequence.

As a simple example, let $p_n = 1-\e^{-2^{-n}}$ for all $n\in\N$.  Since
$p_n\sim 2^{-n}$, the series $\sum_n p_n$ is convergent.  In this case, it is
easy to compute the distribution of the index of the last non-zero $X_n$.  Set
$M: = \max\{\, n\in\N \mid X_n \ne 0 \,\}$, then we have
\begin{align*}
  \PP(M\le m)
  = \prod_{n=m+1}^\infty (1-p_n)
  &=\prod_{n=m+1}^\infty \e^{-2^{-n}} \\
  &=\exp\left(-2^{-m}\sum_{n=1}^\infty 2^{-n} \right)
  = \e^{-2^{-m}}.
\end{align*}
Note that for $m=0$ (which means that all $X_n$ are zero) this gives a
correct answer, $\PP(M=0)=\e^{-1}$. More generally, for $m\ge 1$
\begin{equation}\label{eq:M}
  \PP(M=m)
  = \e^{-2^{-m}}-\e^{-2^{-(m-1)}}
  = \e^{-2^{-m}}\left(1-\e^{-2^{-m}}\right).
\end{equation}
Hence, we can first sample the random number $M=m$ according to its
distribution~\eqref{eq:M} and then sample independent Bernoulli random
variables $X_1,\dots,X_{m-1}$ (with $X_m = 1$ by definition of $M$). This gives
an exact sampling algorithm, but to generalise it we need to be able to
calculate explicitly the infinite products of the form $\prod_{n>m} (1-p_n)$,
which is not always possible.

As a pragmatic alternative, we consider a straightforward but approximate algorithm.

\medbreak

\noindent
\textbf{Algorithm~0.}
\begin{enumerate}[itemsep=.5ex,parsep=0pt,topsep=.5ex]
\item Choose $N$ large.
\item For $n \in \{1, \ldots, N\}$, generate $X_n \in \{0, 1\}$
  independently with $\PP(X_n = 1) = p_n$.
\item Set $X_n = 0$ for all $n > N$.
\end{enumerate}

\medbreak

In the output of Algorithm~0, the events $X_n$ are independent.
For $n \leq N$ we have $\PP(X_n = 1) = p_n$ as required, but
for $n > N$ we have $\PP(X_n = 1) = 0$.  If we choose $N$ large enough that
$p_n$ is `small' for all $n > N$, we can use Algorithm~0 as an approximate
solution to our problem.
Our aim is to present an algorithm (using only finitely many steps)
which can
exactly sample from the full sequence of events, corresponding to the case $N=\infty$
in Algorithm~0.

Sequences of this kind appear in the study of random combinatorial
structures.  An integer partition can be described by the numbers of
parts of each size and, for suitable probability measures on the set of
all partitions, these numbers are independent random variables of which
only finitely many are non-zero
\citep{Fristedt1993,Vershik1996,Bogachev2015}.  Convex lattice polygonal
lines can be described in the same way, using the multiplicities of the
sides in each possible direction \citep{Sinai1994,BogZar2011}.  A further
example is the occupancy problem with infinitely many boxes, where a
Poisson distributed number of balls makes the numbers of balls in
different boxes independent \citep{GneHanPit2007,BogGneYak2008}.  Our
algorithm is based on this last observation.

For combinatorial structures, exact sampling is often considered in the
framework of Boltzmann samplers, introduced by \citet{duchon_2004}.  Such
a sampler generates an object of random size, in such a way that all
objects of the same size are equally likely.  \citet{flajolet_2007} extend
the framework to unlabelled structures, including multisets and powersets;
their samplers assume that an oracle for the values of the generating
functions is available.
If an object of a given size is required, the output of the sampler can be
rejected until the size is correct, and \citet{ArrDes2016} and
\citet{desalvo_2017a} show how the cost of this rejection step can be
reduced.  Our method uses two older ideas.  The first idea is the thinning
method of \citet{lewis_1979}, which simulates a Poisson process by
generating the points of a process with larger intensity and then
deleting some of these points at random.  The second idea is retrospective
sampling, which \citet{BesPapRob2006} use for the exact simulation of
a class of diffusions: there, an infinite-dimensional object is simulated
exactly by generating only the finitely many random variables which are
needed to determine the output.

A related construction appears in parallel work by
\citet{peyen_2026}, who considers the probabilities $z^n / (1+z^n)$, with
several independent events for each~$n$, in order to obtain Boltzmann
samplers for powersets.  There, the independence of the events is derived
from an occupancy model with a Poisson distributed number of balls, and
the rejection step is justified using the thinning method for Poisson
processes.  \citet{peyen_2026} also shows that, in such an occupancy
model with infinitely many boxes, the Poisson distribution is the only distribution for the number
of balls which makes the events uncorrelated.  The focus of the present article is different:
we state the method for arbitrary summable probabilities~$p_n$, give a
self-contained proof based on Poisson processes, and consider examples
outside the setting of Boltzmann samplers.

Section~\ref{sec:exact} introduces the main results of the article,
namely Algorithm~1, the correctness proof in Lemma~\ref{L:correctness}
and the efficiency analysis in Lemma~\ref{L:complexity}.
In section~\ref{sec:applications} we apply the
method to three examples: the prime factorisation of zeta-distributed
random integers, random integer partitions and convex lattice polygonal
lines.  In the first of these examples the probabilities $p_n$ decay so
slowly that Algorithm~0 is not feasible when the parameter of the
distribution is close to~$1$, whereas the expected cost of the exact
algorithm is polynomial.  To our knowledge, this gives the first exact
sampler for factored, zeta-distributed integers.

\section{Exact Sampling}
\label{sec:exact}

The following algorithm gives an exact solution for the sampling problem
described in the introduction.

\subsection{The Algorithm}

\medbreak

\textbf{Algorithm~1.}
Choose probability weights $(\pi_n)_{n\in\N}$ such that $\pi_n > 0$ for all
$n\in\N$ and such that there exists an $a\in\R$ with $a_n := -\log(1-p_n)/\pi_n
\leq a$ for all $n\in\N$.  Then perform the following steps:
\begin{enumerate}[itemsep=1ex,parsep=0pt,topsep=1ex]
\item Generate $K \sim \Poiss(a)$.
\item Generate $N_1, \ldots, N_K \in\N$, independently, such that
  $\PP(N_k = n) = \pi_n$ for all $n\in\N$ and all $k \in \{1, \ldots, K\}$.
\item Generate $U_1, \ldots, U_K \in [0, a]$, independently, uniformly
  distributed.
\item For all $n\in\N$, set $X_n = 1$ if there is a
  $k \in \{1, \ldots, K\}$ with $N_k = n$ and $U_k \leq a_{n}$, and
  set $X_n = 0$ otherwise.
\end{enumerate}

\medbreak

\begin{figure}[t]
  \centering
  \includegraphics{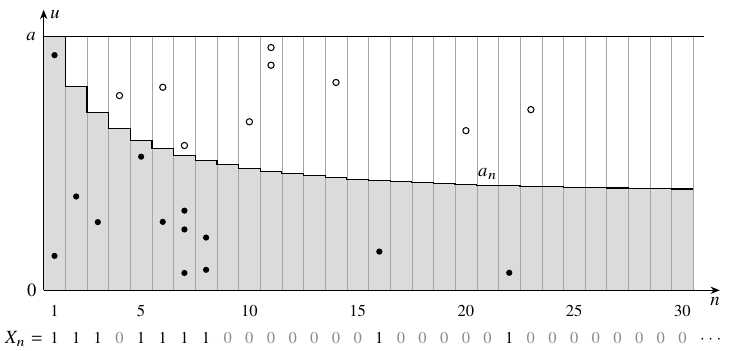}
  \caption{\label{fig:strip}One run of Algorithm~1 for $p_n = z^n$ and
    $\pi_n = z^{n-1}(1-z)$ with $z = 0.9$, as in Example~\ref{ex:geom}.
    Each point $(N_k, U_k)$ is shown as a circle at height~$U_k$ in
    column~$N_k$.  The circle is filled if $U_k \leq a_{N_k}$, and we have
    $X_n = 1$ if and only if column~$n$ contains at least one filled
    circle.}
\end{figure}

Figure~\ref{fig:strip} illustrates the algorithm.  Steps 1 to~3 generate
random points $(N_k, U_k)$ in the strip $\N \times [0, a]$, and step~4
accepts the points which satisfy $U_k \leq a_{N_k}$.  We have $X_n = 1$
if and only if at least one of the accepted points has first
coordinate~$n$.

To use the algorithm, we need to choose probability weights $\pi =
(\pi_n)_{n\in\N}$ which satisfy the assumptions stated at the start of the
algorithm.  Under the assumptions $p_n \in [0, 1)$ and $\sum_{n\in\N} p_n <
\infty$, a suitable distribution~$\pi$ always exists. For example, it is easy
to show that
\begin{equation*}
  a^* := - \sum_{n\in\N} \log(1-p_n) < \infty
\end{equation*}
and thus $\pi_n = -\log(1-p_n) / a^*$ is a valid choice.

The following lemma shows that the algorithm works as intended.

\begin{lemma}\label{L:correctness}
  Given $(p_n)_{n\in\N}$ with $p_n \in [0, 1)$ for all $n\in\N$ and
  $\sum_{n\in\N} p_n < \infty$,
  let $(X_n)_{n\in\N}$ be the sequence of random variables
  constructed by Algorithm~1.
  Then the $X_n$ are independent, with $\PP(X_n = 1) = p_n$ for all~$n\in\N$.
\end{lemma}

The proof, loosely inspired by \citet{BesPapRob2006}, is a nice application of
Poisson processes.  For reference we summarise the required results here,
leaving measurability considerations aside. Full details can be found, for
example, in the monograph by \citet[sections 2.1 and~2.5]{Kingman93}. In the
following, $\#M$ denotes the number of elements in a finite or countable
set~$M$.

\begin{definition}
  Let $S$ be a set and $\Lambda$ a measure on $S$ such that
  $\Lambda(\{x\}) = 0$ for all $x\in S$.
  A {\em Poisson process} on $S$ with intensity measure $\Lambda$
  is a random subset $\Pi$ of $S$ such that
  \begin{enumerate}
  \item for each $m\in\N$ and any disjoint sets $A_1, \ldots, A_m \subseteq S$,
    the random variables $\#(A_1\cap\Pi), \ldots, \#(A_m\cap \Pi)$ are
    independent, and
  \item for every $A \subseteq S$, the random variable $\#(A\cap \Pi)$
    is Poisson distributed with parameter~$\Lambda(A)$.
  \end{enumerate}
\end{definition}

In the proof of Lemma~\ref{L:correctness}, below, we will use the fact
that the $Z_k := (N_k, U_k)$ constructed in the algorithm form a Poisson process.
This follows from the same argument as is used in the
existence proof of Poisson processes: Assume that
$\mu := \Lambda(S) < \infty$
and that $K, Z_1, Z_2, \ldots$ are independent random variables such that
$K \sim \Poiss(\mu)$
and
\begin{equation*}
  P\bigl( Z_k \in A \bigr)
  = \frac{\Lambda(A)}{\Lambda(S)}
\end{equation*}
for all $A \subseteq S$.
Define $\Pi := \bigl\{ Z_k \bigm| k = 1, \ldots, K \bigr\}$.
Let $A_1, \ldots, A_m \subseteq S$ be disjoint.  Then
\begin{align*}
  \hskip5mm&\hskip-5mm
    P\bigl(\#(\Pi\cap A_1) = k_1, \ldots, \#(\Pi\cap A_m) = k_m \bigr) \\
  &= \sum_{k = \sum k_i}^\infty P\bigl(\#(\Pi\cap A_1) = k_1, \ldots, \#(\Pi\cap A_m) = k_m \bigm| K=k \bigr) \cdot P( K = k) \\
  &= \cdots \\
  &= \prod_{i=1}^m \e^{-\Lambda(A_i)} \frac{\Lambda(A_i)^{k_i}}{k_i!}.
\end{align*}
(We leave some of the steps as an exercise to the reader.)
This shows that the $\#(\Pi\cap A_i)$ are independent and Poisson distributed,
and thus that $\Pi$ is indeed a Poisson process.

\begin{proof}[Proof of Lemma~\ref{L:correctness}]
  Let $S := \N \times [0, a]$ and consider the product measure
  $\Lambda = \pi \times \Leb$ on~$S$, where $\Leb$
  denotes the Lebesgue measure on $[0, a]$.  Then we have
  $\Lambda(S) = a < \infty$ and $\Lambda$ has no atoms.
  For the purpose of the proof we extend the $U_i$ and $N_i$
  from the algorithm
  to infinite sequences $(U_i)_{i\in\N}$ and $(N_i)_{i\in\N}$,
  where all $N_i, U_i$ are independent of each other and of $K$.
  From the algorithm we recall that
  $\PP(N_i = n) = \pi_n$ and that the $U_i$ are uniformly distributed on $[0, a]$.
  Since
  \begin{equation*}
    P\bigl( (N_i, U_i) \in A \bigr)
    = \sum_{n=1}^\infty \pi_n \frac{\Leb( \{\, u \mid (n,u) \in A \,\})}{a}
    = \frac{\Lambda(A)}{\Lambda(S)}
  \end{equation*}
  for all $i\in\N$ and all $A\subseteq S$,
  the construction preceding this proof
  shows that the random set
  $\Pi := \bigcup_{k=1}^K \bigl\{ (N_k, U_k) \bigr\} \subset S$
  is a Poisson process with intensity measure~$\Lambda$ on~$S$.

  By definition of $a_n$ we have $a_n = -\log(1-p_n)/\pi_n$.
  Consider the sets
  \begin{equation*}
    S_n
    := \bigl\{ (n, u) \bigm| u \leq a_n \bigr\}
    \subseteq S
  \end{equation*}
  for all $n\in\N$.  In Figure~\ref{fig:strip}, these sets $S_n$ form the
  shaded parts of the columns.  By step~4 of the algorithm we have
  \begin{equation*}
    \{ X_n = 1 \}
    = \{ \Pi\cap S_n \neq \varnothing \}.
  \end{equation*}
  Since $\Pi$ is a Poisson process and since the sets $S_n$
  are disjoint, it follows that the random variables $X_n$ are independent.

  To show that the $X_n$ have the correct distribution, we first write
  \begin{equation*}
    \lambda_n
    := \Lambda(S_n)
    = \pi_n \cdot \mathrm{Leb}\bigl( [0, a_n] \bigr)
    = \pi_n a_n
    = -\log(1 - p_n).
  \end{equation*}
  Using the properties of a Poisson process again we get
  \begin{equation*}
    \#\bigl( \Pi\cap S_n \bigr)
    \sim \mathrm{Poiss}\bigl( \Lambda(S_n) \bigr)
    = \mathrm{Poiss}\bigl( \lambda_n \bigr)
  \end{equation*}
  and thus
  \begin{align*}
    \PP(X_n = 1)
    &= \PP\bigl( \Pi\cap S_n \neq \varnothing \bigr) \\
    &= 1 - \PP\bigl( \# (\Pi\cap S_n) = 0 \bigr) \\
    &= 1 - \exp(-\lambda_n) \\
    &= p_n,
 \end{align*}
 as required. This completes the proof.
\end{proof}

The following two examples illustrate the algorithm.  In the first
example every point is accepted, while in the second example
some of the points are rejected.

\begin{example}\label{ex:intro}
  In the introduction we considered the case
  $p_n = 1 - \e^{-2^{-n}}$ for all $n\in\N$.  Here we have
  $-\log(1-p_n) = 2^{-n}$ and thus $a^* = 1$.  We choose
  $\pi_n = -\log(1-p_n)/a^* = 2^{-n}$, so that $\pi$ is the geometric
  distribution with parameter~$1/2$.  This gives $a_n = 1$ for all
  $n\in\N$, and we can choose~$a = 1$.  Since the condition
  $U_k \leq a_{N_k}$ in step~4 is then always satisfied, the $U_k$ are
  not needed and the algorithm takes the following form: we generate
  $K \sim \Poiss(1)$, then we generate $N_1, \ldots, N_K$ independently
  from the geometric distribution~$\pi$, and finally we set $X_n = 1$ if
  and only if $n \in \{N_1, \ldots, N_K\}$.
\end{example}

\begin{example}\label{ex:geom}
  Let $z\in (0,1)$ and $p_n = z^n$.
  Then we have $\sum_{n\in\N} p_n = z/(1-z) < \infty$, so the method can be applied.
  Choose $\pi$ to be the geometric
  distribution with parameter $1-z$, that is, $\pi_n = z^{n-1} (1-z)$
  for all $n\in\N$.  Then we have
  \begin{equation*}
    a_n
    = - \frac{\log(1-p_n)}{\pi_n}
    = - \frac{\log\bigl( 1 - z^n \bigr)}{z^{n-1} (1-z)}
    = \frac{z}{1-z} \, f(z^n),
  \end{equation*}
  where $f(y) = -\log(1-y)/y$ for all $y\in(0,1)$.  Since
  $f(y) = \sum_{j=0}^\infty y^j/(j+1)$, the function $f$ is monotonically
  increasing.  Thus, the sequence $(a_n)_{n\in\N}$ is monotonically
  decreasing and we can choose
  \begin{equation}\label{eq:a-for-geom}
    a
    := a_1
    = - \frac{\log\bigl( 1 - z \bigr)}{1-z}.
  \end{equation}
  Using this value of~$a$, it is easy to implement the method on a computer.
  One run of the method, for $z = 0.9$, is shown in Figure~\ref{fig:strip}.
\end{example}

\subsection{Efficiency}
\label{sec:efficiency}

When
implementing the algorithm on a computer, the proposal distribution~$\pi$ should be chosen so
that the samples $N_1, \ldots, N_K \sim \pi$ in step~2 can be
generated efficiently.

\begin{lemma}\label{L:complexity}
A run of Algorithm~1 consumes $2K+1$ random numbers, where
$K \sim \Poiss(a)$.  In particular, the expected number of random numbers
consumed is~$2a+1$.
\end{lemma}

\begin{proof}
  Step~1 of the algorithm consumes one random number, namely~$K$, and
  steps 2 and~3 consume $K$ random numbers each.  Step~4 consumes no
  random numbers.  Finally, we have $\E(K) = a$.  This completes
  the proof.
\end{proof}

If each of the random numbers can be generated in a bounded amount of time,
and if the same holds for the evaluation of~$a_n$, the lemma shows that the
expected running time and the expected memory use of Algorithm~1 are bounded
by a constant times~$1+a$.
The algorithm works for any value of the upper bound~$a$, but in the light of Lemma~\ref{L:complexity}
the bound $a$ should be chosen as small as
possible.

The constant $a^*$, defined before Lemma~\ref{L:correctness}, is a
lower bound for~$a$.  Since $\pi_n a_n = -\log(1-p_n)$ and $a_n \leq a$
for all $n\in\N$, we have
\begin{equation*}
  a^*
  = \sum_{n\in\N} \pi_n a_n
  \leq a \sum_{n\in\N} \pi_n
  = a
\end{equation*}
for every choice of $\pi$ and~$a$ which satisfies the conditions of
Algorithm~1.  Equality holds if and only if $a_n = a$ for all $n\in\N$,
and this is the case only for the choice $\pi_n = -\log(1-p_n)/a^*$.
In the general case we have
\begin{equation*}
  \PP\bigl( U_k \leq a_{N_k} \bigr)
  = \sum_{n\in\N} \pi_n \frac{a_n}{a}
  = \frac{a^*}{a}
\end{equation*}
for every~$k$,
and thus the ratio $a^*/a$ is the probability that a point is accepted.
Thus, the expected number of accepted points is $\E(K) \, a^*/a = a^*$.
The algorithm thus generates $a$ points on
average in order to obtain $a^*$ accepted points.

The choice $\pi_n = -\log(1-p_n)/a^*$ minimises the cost, but it can only
be used when samples from this distribution can be generated efficiently.
This is the case in Example~\ref{ex:intro}, where we have $a = a^* = 1$
and every point is accepted.  In the other examples in this article we
choose a distribution~$\pi$ which is easy to sample from, and we then
have~$a > a^*$.

\begin{example}
To compare efficiency of the algorithms 0 and~1, we
use the number of random samples
consumed as a measure for the resource cost of the methods here.
Clearly, Algorithm~0 consumes exactly $N$ random variables, and the value~$N$
increases as the probability of missing one or more events decreases.
From Lemma~\ref{L:complexity} we know that Algorithm~1
uses $2K+1$ samples where $\E(K) = a$, so
the expected number of random variables used by Algorithm~1 is $2a+1$.  Algorithm~1 is exact.

For the situation of Example~\ref{ex:geom},
with $z = 0.9$, equation~\eqref{eq:a-for-geom}
gives $a = 23.026$ and we get the numerical values
\begin{center}
\begin{tabular}{lcc}
method & cost & $\PP(\mathrm{error})$ \\ \hline
Algorithm~0, $N = 50$ & $50$ & $4.538\%$ \\
Algorithm~0, $N = 65$ & $65$ & $0.951\%$ \\
Algorithm~0, $N = 87$ & $87$ & $0.094\%$ \\[1ex]
Algorithm~1 & $\hskip-12pt\approx47$ & exact
\end{tabular}
\end{center}
Thus, for the example the exact method will be faster than any useful version
of Algorithm~0.
\end{example}

\section{Applications}
\label{sec:applications}

In this section we apply Algorithm~1 to three examples.

\subsection{Random Factored Integers}
\label{sec:integers}

For $s > 1$, a random integer $Z$ is zeta-distributed, if it satisfies
\begin{equation*}
  \PP( Z = n )
  = n^{-s} / \zeta(s)
\end{equation*}
for all $n\in \N$,
where $\zeta(s) = \sum_{n=1}^\infty n^{-s}$ is the Riemann zeta function.
We can construct zeta-distributed values by sampling the prime factorisation of~$Z$
as follows: for every prime $q\in\N$ let $\nu_q$ be geometrically distributed with
\begin{equation*}
  \PP( \nu_q = k )
  = q^{-sk} (1 - q^{-s})
\end{equation*}
for all $k\in \N_0$, independently for different~$q$.  Let
\begin{equation*}
  Z := \prod_q q^{\nu_q},
\end{equation*}
where the product is finite since, with probability one, only
finitely many $\nu_q$ are non-zero.
Then, if $n$ has prime factorisation $n = \prod_q q^{\tilde\nu_q}$, we find
\begin{equation*}
  \PP(Z = n)
  = \prod_q \PP\bigl( \nu_q = \tilde\nu_q \bigr)
  = \prod_q q^{-s\tilde\nu_q} (1 - q^{-s})
  = n^{-s} \prod_q (1 - q^{-s}).
\end{equation*}
Using the Euler product formula we get
\begin{equation*}
  \PP(Z = n)
  = n^{-s} \frac{1}{\zeta(s)}
\end{equation*}
and thus $Z$ is zeta-distributed.
This representation is classical: \citet{golomb_1970} shows that, for a
zeta-distributed integer, divisibility by different primes is independent
and that the exponent of each prime is geometrically distributed.
A recent use of the representation is the `zeta process' of
\citet{alexeev_2026}, which couples the $\nu_q$ for different values of~$s$.
Our aim here is to sample $Z$ together with its prime factorisation,
\textit{i.e.}\ to sample the~$\nu_q$.

Sampling $Z$ alone, without its factorisation, is a standard problem,
and \citet[section~X.6.1]{devroye_1986} gives a rejection method for the
zeta distribution, which is called the Zipf distribution there.
Recovering the factorisation afterwards requires factoring~$Z$, and this
is not feasible when $s$ is close to~$1$, because $Z$ is then typically
very large.  Indeed, we have
\begin{equation*}
  \E(\log Z)
  = \frac{1}{\zeta(s)} \sum_{n=1}^\infty n^{-s} \log n
  = -\frac{\zeta'(s)}{\zeta(s)},
\end{equation*}
and this expression behaves like $1/(s-1)$ as $s \downarrow 1$.

Alternatively, we could use Algorithm~0, sampling $\nu_q$ for all primes
$q \leq N$ and setting $\nu_q = 0$ for all primes $q > N$.  The result is
wrong whenever a prime larger than $N$ divides~$Z$, and we have
\begin{equation*}
  \PP\bigl( \nu_q > 0 \mbox{ for some } q > N \bigr)
  = 1 - \prod_{q > N} (1 - q^{-s})
  \approx \sum_{q > N} q^{-s}.
\end{equation*}
By the prime number theorem, the sum on the right-hand side decays only
like $N^{-(s-1)} / \log N$ as $N\to\infty$.  In contrast to
Example~\ref{ex:geom}, where the corresponding probabilities decay
geometrically, a small error probability thus requires a very large
value of~$N$ when $s$ is close to~$1$.

The corresponding problem for uniformly distributed integers was solved by
\citet{bach_1988} and, with a simpler algorithm, by \citet{kalai_2003}.
Both methods generate a random integer together with its prime
factorisation, using primality tests in place of factorisation.  \citet{bodini_2012} construct
samplers for distributions which are defined by Dirichlet series.  Their
samplers cover the zeta distribution without the factorisation as well as
the distribution of a single random prime, but the multiset construction
which would be needed for the prime factorisation is deferred to an
extended version of their article.
To our knowledge, this extended version has not yet appeared.

To apply Algorithm~1, we note that a prime $q$ divides $Z$ if and only if
$\nu_q > 0$.  These events are independent, with
\begin{equation*}
  \PP(\nu_q > 0)
  = q^{-s}
\end{equation*}
for all primes~$q$, and we have $\sum_q q^{-s} \leq \zeta(s) - 1 < \infty$.
Thus, the set of primes which divide~$Z$ is a random set of the kind
considered in this article.  Once this set is known, the exponents are easy
to generate, since
\begin{equation}\label{eq:zeta-exponent}
  \PP( \nu_q = k \,|\, \nu_q > 0)
  = q^{-s(k-1)} (1 - q^{-s})
\end{equation}
for all $k\in\N$ is a geometric distribution.

The index set of the events is now the set of primes instead of~$\N$, and
there is no convenient probability distribution $\pi$ on the primes which
we could use in Algorithm~1.  We solve this problem by extending the index
set: we set $p_n = n^{-s}$ for all integers $n\geq 2$, whether prime or
not, apply Algorithm~1 to this larger family of events and, at the end,
discard all indices which are not prime.  As a minor optimisation,
we treat the prime~$2$ separately, by sampling $\nu_2$
directly, and then restrict the algorithm to the odd integers $n\geq 3$.

To choose the proposal distribution~$\pi$, let $r = 2^{-(s-1)}$ and, for $k\in\N$, let $B_k$ be the
set of all odd integers $n$ with $2^k < n < 2^{k+1}$.  The sets $B_k$ are
disjoint, each $B_k$ has $2^{k-1}$ elements, and every odd integer
$n\geq 3$ is contained in exactly one of the~$B_k$.  We define
\begin{equation*}
  \pi_n
  = (1-r) \, r^{k-1} \, 2^{-(k-1)}
  = (1-r) \, 2^{-s(k-1)}
\end{equation*}
for all $n\in B_k$ and all $k\in\N$.  Since the weights of the elements
of $B_k$ sum to $(1-r) r^{k-1}$, this defines a probability distribution.
To sample from $\pi$, we can first generate $k\in\N$ with probabilities
$(1-r) r^{k-1}$ and then choose $n$ uniformly from~$B_k$, which requires
$k-1$ random bits.

To find the constant~$a$, we again use the fact that the function
$f(y) = -\log(1-y)/y$ from Example~\ref{ex:geom} is monotonically
increasing on~$(0,1)$.  For $n\in B_k$ we have $n > 2^k$ and $n \geq 3$,
and thus
\begin{equation*}
  a_n
  = \frac{-\log(1-n^{-s})}{\pi_n}
  = f(n^{-s}) \, \frac{n^{-s} \, 2^{s(k-1)}}{1-r}
  \leq \frac{f(3^{-s}) \, 2^{-s}}{1-r}
  =: a.
\end{equation*}
As $s\downarrow 1$ we have $1-r \sim (s-1)\log 2$, and thus $a$ behaves
like $c/(s-1)$ where $c = 3\log(3/2) / (2\log 2) \approx 0.88$.
Combining these choices with Algorithm~1, we obtain the following
method, which we call the \emph{zeta sampler}.

\medbreak

\noindent
\textbf{Zeta sampler.}
Let $r$, $B_k$, $a_n$ and $a$ be as above.
\begin{enumerate}[itemsep=1ex,parsep=0pt,topsep=1ex]
\item Generate $\nu_2 \in \N_0$ with
  $\PP(\nu_2 = k) = 2^{-sk} (1 - 2^{-s})$ for all $k\in\N_0$.
\item Generate $K \sim \Poiss(a)$.
\item For $i \in \{1, \ldots, K\}$, independently, generate $k_i\in\N$ with
  $\PP(k_i = k) = (1-r) r^{k-1}$ for all $k\in\N$, choose $N_i$ uniformly
  from $B_{k_i}$, and generate $U_i$ uniformly distributed on~$[0, a]$.
\item Let $I$ be the set of all $N_i$ such that $U_i \leq a_{N_i}$ and
  $N_i$ is prime.
\item For every $q\in I$, generate $\nu_q \in \N$ with the
  probabilities given in~\eqref{eq:zeta-exponent}.  Set $\nu_q = 0$ for all odd
  primes $q \notin I$.
\end{enumerate}

The output of the zeta sampler is the family $(\nu_q)$ of exponents, and
thus the prime factorisation of $Z = \prod_q q^{\nu_q}$.

\medbreak

In step~4, the condition $U_i \leq a_{N_i}$ should be checked first, so
that the more expensive primality test is only needed for those candidates
which satisfy this condition.  By Lemma~\ref{L:complexity}, the expected
number of candidates $N_i$ is~$a$.  A candidate equals $n$ and satisfies
$U_i \leq a_n$ with probability $\pi_n a_n / a$, and thus the expected
number of primality tests is
\begin{equation*}
  \sum_{\substack{n \geq 3 \\ n \text{ odd}}} -\log(1 - n^{-s})
  \approx \sum_{\substack{n \geq 3 \\ n \text{ odd}}} n^{-s}
  = (1 - 2^{-s}) \zeta(s) - 1,
\end{equation*}
which behaves like $1/\bigl(2(s-1)\bigr)$ as $s\downarrow 1$.
The sum on the left-hand side is the constant $a^*$ from
section~\ref{sec:efficiency} for the events with odd index $n\geq 3$.
Since $a$ behaves like $c/(s-1)$, the probability $a^*/a$ that a candidate
is accepted converges to
$1/(2c) = \log 2 / \bigl( 3\log(3/2) \bigr) \approx 0.57$ as
$s\downarrow 1$.  A candidate
from $B_k$ has $k+1$ binary digits, and the expected number of digits is
$1 + 1/(1-r) \approx 1.44 / (s-1)$.  Since the time needed to test an
integer for primality is polynomial in the number of
digits~\citep{agrawal_2004}, the expected running time of the zeta sampler is
bounded by a polynomial in~$1/(s-1)$.  This matches the size of the
output, since $\E(\log Z)$ is also of order~$1/(s-1)$.

For comparison, we consider the cut-off $N$ which Algorithm~0 needs in
order to achieve an error probability of~$1\%$.  By the prime number
theorem we have
\begin{equation*}
  \sum_{q > N} q^{-s}
  \approx \int_M^\infty \frac{x^{-s}}{\log x} \, \rmd x
  = E_1\bigl( (s-1) \log N \bigr),
\end{equation*}
where $E_1(x) = \int_x^\infty \e^{-t} / t \, \rmd t$ is the exponential
integral.  The error probability equals $1\%$ when the sum equals
$-\log(0.99)$, and since $E_1(3.207) \approx -\log(0.99)$, the required
cut-off is $N \approx \exp\bigl( 3.207 / (s-1) \bigr)$, which grows
exponentially in~$1/(s-1)$.  We get the following numerical values.
\begin{center}
\begin{tabular}{lcccc}
  & & expected number of & expected number of & cut-off $N$ \\
$s$ & $a$ & primality tests & binary digits & for Algorithm~0 \\ \hline
$1.5$ & $1.34$ & $0.72$ & $4.4$ & $601$ \\
$1.1$ & $8.27$ & $4.74$ & $15.9$ & $10^{14}$ \\
$1.01$ & $87.2$ & $49.8$ & $145.8$ & $10^{139}$ \\
$1.001$ & $876.9$ & $499.8$ & $1444.2$ & $10^{1393}$
\end{tabular}
\end{center}
The cut-off for $s = 1.5$ was determined by direct computation, the other
values of $N$ were obtained using the approximation above.  For $s = 1.01$,
for example, the zeta sampler tests on average $50$ integers with about $146$
binary digits each for primality, whereas Algorithm~0 would need to consider
all primes up to~$10^{139}$.

We conclude this section with three remarks.  First, the zeta sampler is exact
only if the primality test in step~4 is.  If a probabilistic method like
the test of \citet{rabin_1980} is used instead, each composite candidate is
accepted as a prime factor with a probability of less than~$4^{-k}$, where
$k$ is the number of rounds of the test.  Second, the proposal distribution used here is
not the only possible choice.  For example, we could use the zeta
distribution itself, conditioned on $n\geq 2$, as the proposal for all
primes including~$2$.  This corresponds to
$\pi_n = n^{-s} / \bigl(\zeta(s) - 1\bigr)$ and leads to
$a = \bigl(\zeta(s) - 1\bigr) f(2^{-s}) \approx 1.39/(s-1)$, which is larger
than the value obtained above.
Finally, the same method applies to related distributions.  If we replace
the $\nu_q$ by independent random variables with values in $\{0, 1\}$ and
$\PP(\nu_q = 1) = q^{-s} / (1 + q^{-s})$, then $Z$ is square-free with
$\PP(Z = n)$ proportional to $n^{-s}$ for all square-free~$n$.  Since
$-\log\bigl(1 - \PP(\nu_q = 1)\bigr) = \log(1 + q^{-s}) \leq q^{-s}$, we can
use the same $\pi$ with $a = 2^{-s}/(1-r)$.  Step~1 changes accordingly,
and step~5 is no longer needed.
This is a multiplicative analogue of the powersets considered
by~\citet{peyen_2026}.

\subsection{Random Integer Partitions}
\label{S:int-part}

As another application of Algorithm~1 we consider how to sample
from the uniform distribution on all integer partitions of a given size.
For an integer partition $\lambda_1 + \dots + \lambda_m$,
we define
\begin{equation*}
  \nu_j
  = \sum_{i=1}^\infty 1_{\{ \lambda_i = j \}}.
\end{equation*}
Since the order of the $\lambda_i$ does not matter, each integer partition $\lambda =
(\lambda_1, \dots, \lambda_m)$
corresponds to a unique sequence $\nu = (\nu_j)_{j=1}^\infty$.
Let $z \in (0,1)$.
As shown by \citet{Fristedt1993}, choosing $\nu_j \sim \mathrm{Geom}(1-z^j)$ for all $j\in\N$, independently,
leads to a distribution of $\lambda$ where all partitions of a given size
$|\lambda| = \lambda_1 + \cdots + \lambda_m$
have the
same probability: for a fixed partition $\tilde\lambda$, corresponding to a
sequence $\tilde\nu$ we have
\begin{align*}
  \PP(\lambda = \tilde\lambda)
  &= \PP(\nu = \tilde\nu)
  = \prod_{j=1}^\infty \PP(\nu_j = \tilde\nu_j) \\
  &= \prod_{j=1}^\infty z^{j\tilde\nu_j}(1-z^j)
  = z^{\sum_{j=1}^\infty j\tilde\nu_j} \prod_{j=1}^\infty (1 - z^j)
  = z^{|\tilde\lambda|} \prod_{j=1}^\infty (1 - z^j).
\end{align*}

Since
\begin{equation*}
  \sum_{j=1}^\infty \PP(\nu_j > 0)
  = \sum_{j=1}^\infty z^j
  = \frac{z}{1-z}
  < \infty,
\end{equation*}
we can use Algorithm~1 to exactly simulate the set \begin{equation*}
  \mathcal{I}
  = \bigl\{ j\in\N \bigm| \nu_j > 0 \bigr\},
\end{equation*}
as was demonstrated in Example~\ref{ex:geom}.
Then we have
\begin{equation*}
  \PP( \nu_j = k \,|\, j\in \mathcal{I})
  = \PP( \nu_j = k \,|\, \nu_j > 0)
  = \frac{z^{jk}(1-z^j)}{1-(1-z^j)}
  = z^{j(k-1)}(1-z^j)
\end{equation*}
for all $k \in\N$ and thus for all $j\in \mathcal{I}$ we can simply
generate $\nu_j \sim \mathrm{Geom}(1-z^j)+1$, independently.
Finally, for all $j \notin \mathcal{I}$, we set $\nu_j = 0$.

To generate random partitions from the uniform distribution on all
partitions of an integer~$n$, we can now first generate $\nu$ as
described above, then take the corresponding~$\lambda$, and finally
condition on $|\lambda| = n$ by rejecting all generated $\lambda$ with
$|\lambda| \neq n$.

While this method works for every $z \in (0,1)$, we can tune~$z$
to make the final rejection step more efficient.  The average weight of
a partition generated by the algorithm is
\begin{equation*}
  \E\bigl( |\lambda| \bigr)
  = \E\Bigl( \sum_{j=1}^\infty j \nu_j \Bigr)
  = \sum_{j=1}^\infty j \frac{z^j}{1-z^j}.
\end{equation*}
This expression increases to infinity as $z \uparrow 1$.
Setting $\delta := -\log z$ we find
\begin{equation*}
  \E\bigl( |\lambda| \bigr)
  = \frac{1}{\delta^2} \sum_{j=1}^\infty j\delta \frac{\e^{-j\delta}}{1-\e^{-j\delta}} \,\delta
\end{equation*}
and for $z\approx 1$ we can use Riemann approximation to
get
\begin{equation*}
  \E\bigl( |\lambda| \bigr)
  \approx \frac{1}{\delta^2} \int_0^\infty x \frac{\e^{-x}}{1-\e^{-x}} \,\rmd x
  = \frac{\pi^2}{6\delta^2}
  = \frac{\pi^2}{6\log(z)^2}.
\end{equation*}
Thus, for large~$n$, we can choose
$z = \exp\bigl( -\pi/\sqrt{6n} \bigr)$ to get
$\E\bigl( |\lambda| \bigr) \approx n$.

Similarly, the variance is given by
\begin{align*}
  \Var(|\lambda|)&=\sum_{j=1}^{\infty} j^2 \Var (\nu_j) =\sum_{j=1}^{\infty}j^2 \frac{z^j}{(1-z^j)^2} \\
  &= \frac{1}{\delta^3}\sum_{j=1}^{\infty} (j\delta)^2\frac{\e ^{-j \delta}}{(1-\e^{-j\delta})^2}\delta \\
  &\approx \frac{1}{\delta^3}\int_{0}^{\infty}x^2\frac{\e^{-x}}{(1-\e^{-x})^2} \, \rmd x
  = \frac{\pi^2}{3\delta^3}.
\end{align*}
Since $\delta=-\log z = \pi/\sqrt{6n} \propto n^{-1/2}$, we have
\begin{equation*}
  \Var(|\lambda|)=\frac{\pi^2}{3\delta^3}\propto n^{3/2}.
\end{equation*}

Let $\rho$ denote the expected number of attempts needed to obtain a partition
with $|\lambda|=n$. Since the acceptance probability is of order
$1/\mathrm{sd}(|\lambda|) \propto n^{-3/4}$, we have $\rho \propto n^{3/4}$.
The cost of each attempt is determined by the expected number of Poisson
points. Since $z=\e ^{-\pi/\sqrt{6n}}\approx 1-\frac{\pi}{\sqrt{6n}}$, we have
\begin{equation*}
  \E(K)=a = -\frac{\log(1-z)}{1-z} \approx \frac{-\log\frac{\pi}{\sqrt{6n}}}{\frac{\pi}{\sqrt{6n}}}\propto \sqrt{n}\log n.
\end{equation*}
Thus the overall complexity of the rejection sampling algorithm is
\begin{equation*}
\rho \cdot \E(K) \propto n^{3/4} \cdot n^{1/2}\log n = n^{5/4}\log n.
\end{equation*}
The factor $\log n$ is caused by the geometric proposal distribution
from Example~\ref{ex:geom}.  Using the Riemann approximation again, we
find that the lower bound $a^*$ from section~\ref{sec:efficiency}
satisfies
\begin{equation*}
  a^*
  = -\sum_{j=1}^\infty \log(1-z^j)
  \approx \frac{1}{\delta} \int_0^\infty -\log(1-\e^{-x}) \, \rmd x
  = \frac{\pi^2}{6\delta}
  = \pi \sqrt{n/6},
\end{equation*}
and thus the probability $a^*/a$ that a point is accepted is only of
order~$1/\log n$.  A proposal distribution for which $a$ is of the same
order as $a^*$ would reduce the cost to order~$n^{5/4}$.

For comparison, \citet[section~5.1]{ArrDes2016} represent each geometric
random variable as a sum of independent Poisson random variables and
generate all~$\nu_j$ from a single Poisson process.  This requires of
order $\sqrt{n}$ random numbers per attempt, which is the order
of~$a^*$, so that the rejection method
has cost of order~$n^{5/4}$, without the logarithmic factor, and their
probabilistic divide-and-conquer method reduces the number of attempts
further.  For integer partitions, Algorithm~1 thus does not improve on
the existing methods.  We include the example because it is the simplest
application of the algorithm.

\subsection{Convex Polygonal Lines}
\label{sec:polygons}

A convex polygonal line is a path from the origin~$(0,0)$ to some
endpoint~$(m_1, m_2)$, with vertices in~$\mathbb{Z}^2$ and
consisting of segments with strictly
increasing slope.  As observed by \citet{Sinai1994}, such paths are in one-to-one
correspondence with finitely supported functions~$\nu \colon \mathcal{X} \to \N_0$,
where
\begin{equation*}
  \mathcal{X}
  = \bigl\{ (i_1, i_2) \in \N_0^2 \bigm| \gcd(i_1, i_2) = 1 \bigr\}
\end{equation*}
and~$\nu(x)$ gives the multiplicity of the side in direction~$x$.

For parameters $0 < z_1, z_2 < 1$, consider the family of independent random variables
$(\nu_x)_{x\in\mathcal{X}}$ with geometric distributions, given by
\begin{equation}\label{eq:gcd-dist}
  \PP(\nu_{i_1, i_2} = k)
  = z_1^{k\myp i_1}z_2^{k\myp i_2} \bigl( 1- z_1^{i_1}z_2^{i_2} \bigr)
\end{equation}
for all $(i_1, i_2) \in \mathcal{X}$ and all $k\in\N_0$.
The path corresponding to~$\nu$ has endpoint
$(m_1, m_2) = \bigl( \sum_x \nu_x x_1,\, \sum_x \nu_x x_2 \bigr)$, and all paths
with a given endpoint have the same probability: for a fixed~$\tilde\nu$ with
endpoint~$(\tilde m_1, \tilde m_2)$ we have
\begin{align*}
  \PP(\nu = \tilde\nu)
  &= \prod_{x \in \mathcal{X}} z_1^{\tilde\nu_x \myp x_1} z_2^{\tilde\nu_x \myp x_2} (1 - z_1^{x_1} z_2^{x_2}) \\
  &= z_1^{\tilde m_1} z_2^{\tilde m_2} \prod_{x \in \mathcal{X}} (1 - z_1^{x_1} z_2^{x_2}),
\end{align*}
which does not depend on $\tilde\nu$ except via the endpoint of the path.

Our aim is to sample from the random field $(\nu_x)_{x\in\mathcal{X}}$
using Algorithm~1.
Because there is no simple closed form for the normalizing constant of a
probability distribution supported on $\mathcal{X}$, it is difficult to apply
the algorithm directly to the random field on $\mathcal{X}$.  However, since
the distribution $\pi$ factorizes as a product of geometric distributions on
$\N^2$, we can first extend the random field from $\mathcal{X}$ to $\N^2$,
apply Algorithm~1 there, and then restrict the result back to $\mathcal{X}$,
as we did for the primes in section~\ref{sec:integers}.

As before, we use Algorithm~1 to determine which $\nu_x$ should be
non-zero, and then sample each non-zero $\nu_x$ from the
distribution given in~\eqref{eq:gcd-dist}, conditioned on $\nu_x > 0$.
For this, let
\begin{equation*}
  p_{i_1,i_2}
  := P( \nu_{i_1,i_2} \neq 0)
  = 1 - P(\nu_{i_1,i_2} = 0)
  = 1 - \bigl( 1- z_1^{i_1}z_2^{i_2} \bigr)
  = z_1^{i_1}z_2^{i_2}.
\end{equation*}

Here we choose
\begin{equation*}
  \pi_{i_1, i_2}
  = (1-z_1)z_1^{i_1-1} (1-z_2)z_2^{i_2-1}
\end{equation*}
for all $(i_1, i_2) \in \N^2$.  Since $\pi$ is the product of two
geometric distributions, it is clear that $\pi$ is a probability
distribution.  Furthermore, we have
\begin{align*}
  -\frac{\,\log(1-p_{i_1, i_2})}{\pi_{i_1, i_2}}
  &= -\frac{\log\bigl( 1 - z_1^{i_1}z_2^{i_2} \bigr)}{(1-z_1)z_1^{i_1-1} (1-z_2)z_2^{i_2-1}} \\
  &= \frac{z_1z_2}{(1-z_1)(1-z_2)} \cdot
      \frac{-\log\bigl( 1 - z_1^{i_1}z_2^{i_2} \bigr)}{z_1^{i_1} z_2^{i_2}}.
\end{align*}
Since $y \mapsto -\log(1-y)/y$ is monotonically increasing on $(0,1)$ and
$z_1, z_2 \in (0,1)$ we find
\begin{align*}
  -\frac{\,\log(1-p_{i_1, i_2})}{\pi_{i_1, i_2}}
  &\leq \frac{z_1z_2}{(1-z_1)(1-z_2)} \cdot
      \frac{-\log(1 - z_1 z_2)}{z_1z_2} \\
  &= - \frac{\log(1 - z_1 z_2)}{(1-z_1)(1-z_2)} \\
  &=: a
  < \infty.
\end{align*}

Thus, we can apply Algorithm~1 with index set~$\N^2$ and discard all
points $(i_1, i_2)$ with $\gcd(i_1, i_2) \neq 1$.  By
Lemma~\ref{L:complexity}, the expected number of proposals is~$a$.  The set
$\mathcal{X}$ also contains the two points $(1,0)$ and~$(0,1)$, which are
not elements of~$\N^2$.  For these two points we sample $\nu_{1,0}$ and
$\nu_{0,1}$ directly from the geometric distributions
in~\eqref{eq:gcd-dist}.

To generate convex polygonal lines from the uniform distribution on all
such lines with a given endpoint~$(m_1, m_2)$, we can proceed as for
integer partitions in section~\ref{S:int-part}.  We first generate~$\nu$
as described above, then take the corresponding path, and finally reject
the path if it does not end at~$(m_1, m_2)$.  Since all paths with the
same endpoint have the same probability, the accepted path is uniformly
distributed.  This method works for all $z_1, z_2 \in (0,1)$, but the
final rejection step is more efficient if $z_1$ and~$z_2$ are chosen so
that the expected endpoint of the path is close to~$(m_1, m_2)$.

\vskip2\baselineskip

\textbf{A note on the use of AI.} The original draft,
including Sections \ref{sec:exact}, \ref{S:int-part} and~\ref{sec:polygons},
was written by the authors, before large language models existed.
The application which is now Section~\ref{sec:integers} was suggested by AI
(Anthropic's Claude),
and AI was also used to help with the literature review, to polish
the presentation, and to write the Python script which in turn produced
Figure~\ref{fig:strip}.

\bibliography{references}

\end{document}